# On Contact Numbers of Finite Lattice Sphere Packings of 20-27 Balls

István Szalkai

University of Pannonia, Veszprém, Hungary
szalkai@almos.uni-pannon.hu

March 16, 2016

Abstract

Empirical constructions having maximal contact numbers of unit balls as putative best ones are presented for 20-27 balls.

## 1 Introduction

Consider $n$ balls of the same (unit) radius in the 3-dimensional Euclidean space when $n$ is any but fixed number. If we assume that each pair of balls have at most one common (*tangential*, "*touching*" or "*kissing*") point, the number of such points is called the **contact number** of this ball-configuration. For each fixed number $n \in \mathbb{N}$ we may ask for configurations which have *maximal* contact number, these maximal numbers are denoted by $c(n)$.

Many recent contributions deal with the maximal contact number of balls, using theoretical and empirical approach as well (see eg. *all* the references listed). The possible configurations have many applications e.g. in material science and other fields of applied physics and chemistry.

This enlarged interest lead us to publish our 5 year old empirical investigations. The computer runs investigated all possible cases, so we may think our maximal ones are the real values of $c(n)$. Our investigations are still in progress, so further results will be contributed soon.

Harboth [7] proved that in the plane the hexagonal grid (lattice) is the optimal configuration for congruent circles. This means, that the centers of the circles fit on the grid

$$\left\{ i \cdot \begin{bmatrix} 2 \\ 0 \end{bmatrix} + j \cdot \begin{bmatrix} 1 \\ \sqrt{3} \end{bmatrix} \ : \ i, j \in \mathbb{Z} \right\} \tag{1}$$

For balls in three dimension we suspect that optimal configurations can be found in the 3-dimensional hexagonal grid as follows. First we place some balls in *planes* in the above planar manner, we call this set a **layer**. Second, put such layers onto each other: the layers are translated with the vector $\left[1, \sqrt{1/3}, \sqrt{8/3}\right]^T$ to get the neighbour layer (like the usual placement of melons).

We use the following hexagonal coordinates for the centers of the balls: $P[i, j, k]_H^T$ refers to the $k$-th layer, so $k$ is the vertical coordinate, and $i, j$ are the horizontal ones as in (1). So, $P$ has the the following (usual) Cartesian coordinates (D in the subscript):

for $k$ even :

$$\begin{bmatrix} i \\ j \\ k \end{bmatrix}_H = k \cdot \begin{bmatrix} 0 \\ 0 \\ 2 \end{bmatrix}_D + i \cdot \begin{bmatrix} 2 \\ 0 \\ 0 \end{bmatrix}_D + j \cdot \begin{bmatrix} 1 \\ \sqrt{3} \\ 0 \end{bmatrix}_D \tag{2}$$

for $k = 2l+1$ even :

$$\begin{bmatrix} i \\ j \\ k \end{bmatrix}_H = \begin{bmatrix} 1 \\ \sqrt{1/3} \\ \sqrt{8/3} \end{bmatrix}_D + 2l \cdot \begin{bmatrix} 0 \\ 0 \\ 2 \end{bmatrix}_D + i \cdot \begin{bmatrix} 2 \\ 0 \\ 0 \end{bmatrix}_D + j \cdot \begin{bmatrix} 1 \\ \sqrt{3} \\ 0 \end{bmatrix}_D \tag{3}$$

In the tables in Section 2 we present the centers $\mathbf{F_1}$ through $\mathbf{F_n}$ of the balls of our configurations in hexagonal $[i, j, k]_H^T$ coordinates, which can be transformed to Cartesian ones by (2) and (3) for visualization.

Theorem 4.1 of [5] states for $n \to \infty$ :

$$0.926 < \frac{6n - c(n)}{n^{2/3}} < 7.862$$

The similar (also triangular) *barycentric* coordinate system is also used in chemistry, see e.g. [13].

## 2 Results

First we summarize in a table the maximal contact numbers $c(n)$ (found by computer), which can be considered putative maximal ones:

| $n$ | **20** | **21** | **22** | **23** | **24** | **25** | **26** | **27** |
|---|---|---|---|---|---|---|---|---|
| $c(n)$ | 64 | 67 | 72 | 76 | 80 | 84 | 87 | 90 |
| $c(n)-3n$ | 4 | 4 | 6 | 7 | 8 | 9 | 9 | 9 |
| $\frac{6n - c(n)}{n^{2/3}}$ | 7.6 | 7.7 | 7.6 | 7.6 | 7.7 | 7.7 | 7.8 | 8.0 |

In what follows we give one example for each of the above results. $\mathbf{F_1}$ through $\mathbf{F_n}$ are the centres of the balls, the rows below them are their hexagonal coordinates as defined in the Introduction. Double vertical lines separate the layers, which can be distinguished by the third coordinates. The colors of balls and their connections of separate layers are changing from layer to layer. After the tables all the contacts $(\mathbf{F_i}, \mathbf{F_j})$ of the configuration are listed.

## 20 balls

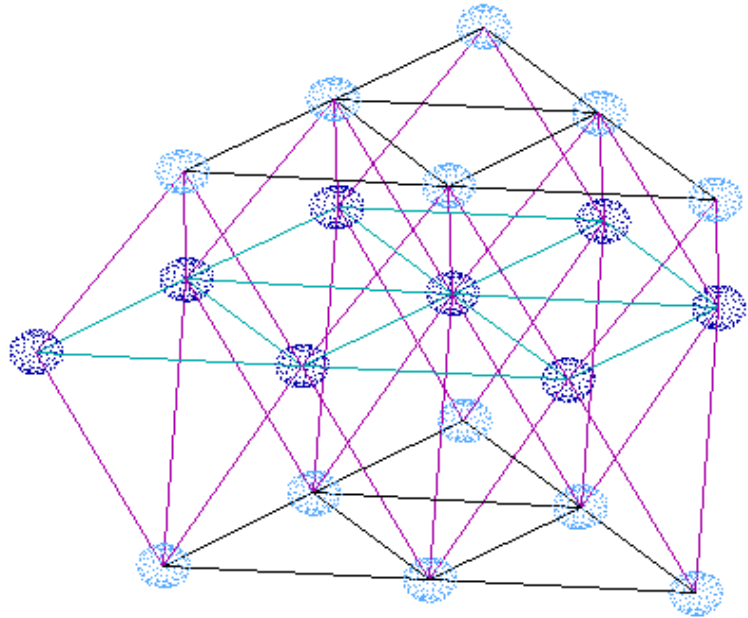

TEX20.acd

| F01 | F02 | F03 | F04 | F05 | F06 | F07 | F08 | F09 | F10 |
|-----|-----|-----|-----|-----|-----|-----|-----|-----|-----|
| 2 | 1 | 2 | 0 | 1 | 2 | 1 | 2 | 0 | 1 |
| 0 | 1 | 1 | 2 | 2 | 2 | 0 | 0 | 1 | 1 |
| 0 | 0 | 0 | 0 | 0 | 0 | 1 | 1 | 1 | 1 |

| F11 | F12 | F13 | F14 | F15 | F16 | F17 | F18 | F19 | F20 |
|-----|-----|-----|-----|-----|-----|-----|-----|-----|-----|
| 2 | 0 | 1 | 2 | 2 | 1 | 2 | 0 | 1 | 2 |
| 1 | 2 | 2 | 2 | 0 | 1 | 1 | 2 | 2 | 2 |
| 1 | 1 | 1 | 1 | 2 | 2 | 2 | 2 | 2 | 2 |

**Contacts:**
(F1,F2), (F1,F3), (F1,F7), (F1,F8), (F2,F3), (F2,F4), (F2,F5), (F2,F7), (F2,F9), (F2,F10), (F3,F5), (F3,F6), (F3,F8), (F3,F10), (F3,F11), (F4,F5), (F4,F9), (F4,F12), (F5,F6), (F5,F10), (F5,F12), (F5,F13), (F6,F11), (F6,F13), (F6,F14), (F7,F8), (F7,F9), (F7,F10), (F7,F15), (F7,F16), (F8,F10), (F8,F11), (F8,F15), (F8,F17), (F9,F10), (F9,F12), (F9,F16), (F9,F18), (F10,F11), (F10,F12), (F10,F13), (F10,F16), (F10,F17), (F10,F19), (F11,F13), (F11,F14), (F11,F17), (F11,F20), (F12,F13), (F12,F18), (F12,F19), (F13,F14), (F13,F19), (F13,F20), (F14,F20), (F15,F16), (F15,F17), (F16,F17), (F16,F18), (F16,F19), (F17,F19), (F17,F20), (F18,F19), (F19,F20),

**Total:** $c(20) = 64$ .

21 balls

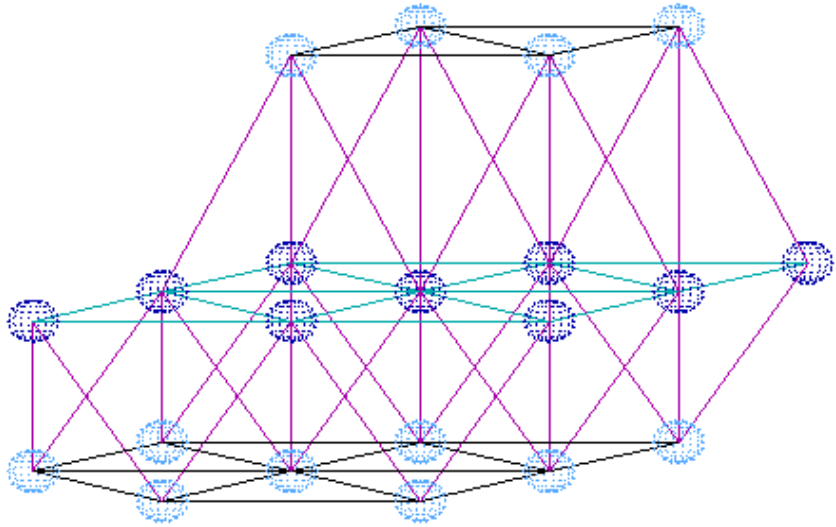

TEX21.acd

| F01 | F02 | F03 | F04 | F05 | F06 | F07 | F08 | F09 | F10 |
|---|---|---|---|---|---|---|---|---|---|
| 1 | 2 | 0 | 1 | 2 | 0 | 1 | 2 | 0 | 1 |
| 0 | 0 | 1 | 1 | 1 | 2 | 2 | 2 | 0 | 0 |
| 0 | 0 | 0 | 0 | 0 | 0 | 0 | 0 | 1 | 1 |

| F11 | F12 | F13 | F14 | F15 | F16 | F17 | F18 | F19 | F20 | F21 |
|---|---|---|---|---|---|---|---|---|---|---|
| 2 | 0 | 1 | 2 | 0 | 1 | 2 | 1 | 2 | 1 | 2 |
| 0 | 1 | 1 | 1 | 2 | 2 | 2 | 1 | 1 | 2 | 2 |
| 1 | 1 | 1 | 1 | 1 | 1 | 1 | 2 | 2 | 2 | 2 |

**Contacts:**
(F1,F2), (F1,F3), (F1,F4), (F1,F9), (F1,F10), (F2,F4), (F2,F5), (F2,F10), (F2,F11), (F3,F4), (F3,F6), (F3,F9), (F3,F12), (F4,F5), (F4,F6), (F4,F7), (F4,F10), (F4,F12), (F4,F13), (F5,F7), (F5,F8), (F5,F11), (F5,F13), (F5,F14), (F6,F7), (F6,F12), (F6,F15), (F7,F8), (F7,F13), (F7,F15), (F7,F16), (F8,F14), (F8,F16), (F8,F17), (F9,F10), (F9,F12), (F10,F11), (F10,F12), (F10,F13), (F10,F18), (F11,F13), (F11,F14), (F11,F19), (F12,F13), (F12,F15), (F12,F18), (F13,F14), (F13,F15), (F13,F16), (F13,F18), (F13,F19), (F13,F20), (F14,F16), (F14,F17), (F14,F19), (F14,F21), (F15,F16), (F15,F20), (F16,F17), (F16,F20), (F16,F21), (F17,F21), (F18,F19), (F18,F20), (F19,F20), (F19,F21), (F20,F21),

**Total:** $c(21) = 67$.

22 balls

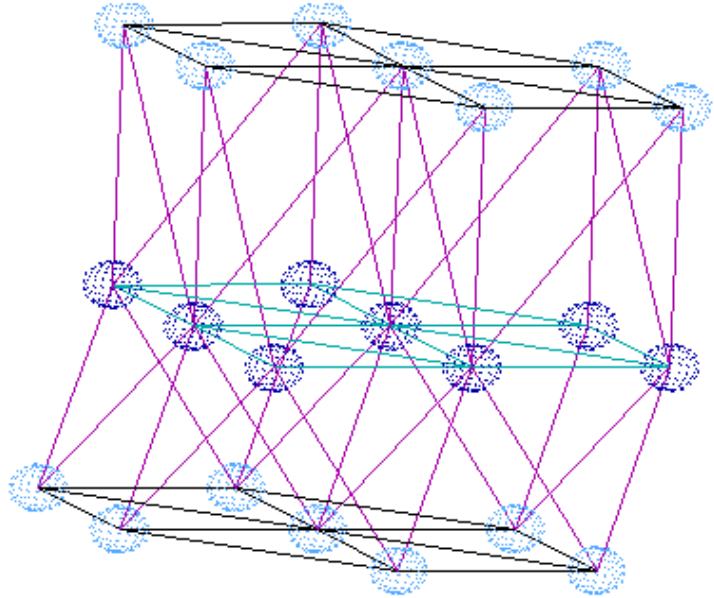

TEX22.acd

| F01 | F02 | F03 | F04 | F05 | F06 | F07 | F08 | F09 | F10 |
|---|---|---|---|---|---|---|---|---|---|
| 1 | 2 | 0 | 1 | 2 | 0 | 1 | 0 | 1 | 2 |
| 0 | 0 | 1 | 1 | 1 | 2 | 2 | 0 | 0 | 0 |
| 0 | 0 | 0 | 0 | 0 | 0 | 0 | 1 | 1 | 1 |

| F11 | F12 | F13 | F14 | F15 | F16 | F17 | F18 | F19 | F20 | F21 | F22 |
|---|---|---|---|---|---|---|---|---|---|---|---|
| 0 | 1 | 2 | 0 | 1 | 1 | 2 | 0 | 1 | 2 | 0 | 1 |
| 1 | 1 | 1 | 2 | 2 | 0 | 0 | 1 | 1 | 1 | 2 | 2 |
| 1 | 1 | 1 | 1 | 1 | 2 | 2 | 2 | 2 | 2 | 2 | 2 |

**Contacts:**

(F1,F2), (F1,F3), (F1,F4), (F1,F8), (F1,F9), (F2,F4), (F2,F5), (F2,F9), (F2,F10), (F3,F4), (F3,F6), (F3,F8), (F3,F11), (F4,F5), (F4,F6), (F4,F7), (F4,F9), (F4,F11), (F4,F12), (F5,F7), (F5,F10), (F5,F12), (F5,F13), (F6,F7), (F6,F11), (F6,F14), (F7,F12), (F7,F14), (F7,F15), (F8,F9), (F8,F11), (F8,F16), (F8,F18), (F9,F10), (F9,F11), (F9,F12), (F9,F16), (F9,F17), (F9,F19), (F10,F12), (F10,F13), (F10,F17), (F10,F20), (F11,F12), (F11,F14), (F11,F18), (F11,F19), (F11,F21), (F12,F13), (F12,F14), (F12,F15), (F12,F19), (F12,F20), (F12,F22), (F13,F15), (F13,F20), (F14,F15), (F14,F21), (F14,F22), (F15,F22), (F16,F17), (F16,F18), (F16,F19), (F17,F19), (F17,F20), (F18,F19), (F18,F21), (F19,F20), (F19,F21), (F19,F22), (F20,F22), (F21,F22),

**Total:** $c(22) = \mathbf{72}$ .

23 balls

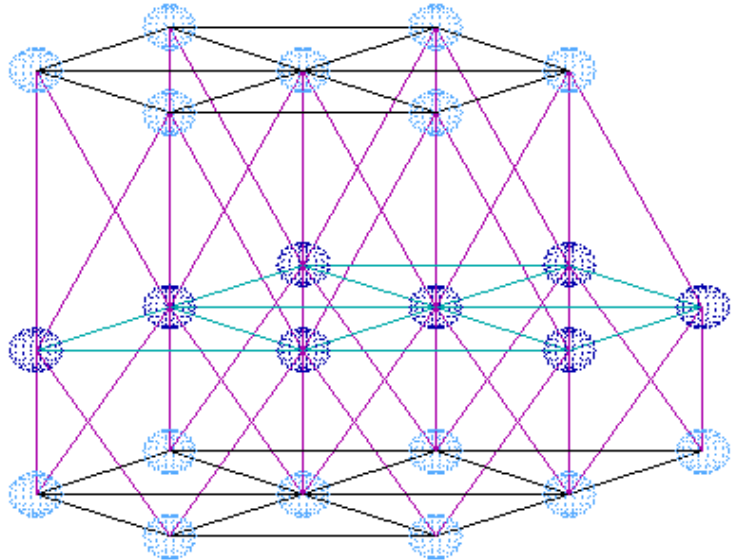

TEX23.acd

| F01 | F02 | F03 | F04 | F05 | F06 | F07 | F08 | F09 | F10 | F11 | F12 |
|---|---|---|---|---|---|---|---|---|---|---|---|
| 1 | 2 | 0 | 1 | 2 | 0 | 1 | 2 | 0 | 1 | 2 | 0 |
| 0 | 0 | 1 | 1 | 1 | 2 | 2 | 2 | 0 | 0 | 0 | 1 |
| 0 | 0 | 0 | 0 | 0 | 0 | 0 | 0 | 1 | 1 | 1 | 1 |

| F13 | F14 | F15 | F16 | F17 | F18 | F19 | F20 | F21 | F22 | F23 |
|---|---|---|---|---|---|---|---|---|---|---|
| 1 | 2 | 0 | 1 | 1 | 2 | 0 | 1 | 2 | 0 | 1 |
| 1 | 1 | 2 | 2 | 0 | 0 | 1 | 1 | 1 | 2 | 2 |
| 1 | 1 | 1 | 1 | 2 | 2 | 2 | 2 | 2 | 2 | 2 |

**Contacts:**
(F1,F2), (F1,F3), (F1,F4), (F1,F9), (F1,F10), (F2,F4), (F2,F5), (F2,F10), (F2,F11), (F3,F4), (F3,F6), (F3,F9), (F3,F12), (F4,F5), (F4,F6), (F4,F7), (F4,F10), (F4,F12), (F4,F13), (F5,F7), (F5,F8), (F5,F11), (F5,F13), (F5,F14), (F6,F7), (F6,F12), (F6,F15), (F7,F8), (F7,F13), (F7,F15), (F7,F16), (F8,F14), (F8,F16), (F9,F10), (F9,F12), (F9,F17), (F9,F19), (F10,F11), (F10,F12), (F10,F13), (F10,F17), (F10,F18), (F10,F20), (F11,F13), (F11,F14), (F11,F18), (F11,F21), (F12,F13), (F12,F15), (F12,F19), (F12,F20), (F12,F22), (F13,F14), (F13,F15), (F13,F16), (F13,F20), (F13,F21), (F13,F23), (F14,F16), (F14,F21), (F15,F16), (F15,F22), (F15,F23), (F16,F23), (F17,F18), (F17,F19), (F17,F20), (F18,F20), (F18,F21), (F19,F20), (F19,F22), (F20,F21), (F20,F22), (F20,F23), (F21,F23), (F22,F23),

**Total:** $c(23) = 76$.

24 balls

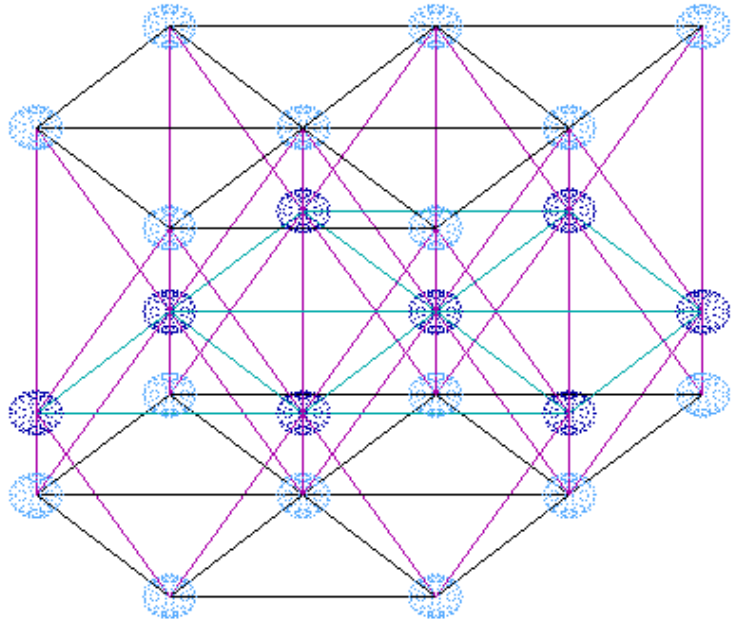

TEX24.acd

| F01 | F02 | F03 | F04 | F05 | F06 | F07 | F08 | F09 | F10 | F11 | F12 |
|-----|-----|-----|-----|-----|-----|-----|-----|-----|-----|-----|-----|
| 1 | 2 | 0 | 1 | 2 | 0 | 1 | 2 | 0 | 1 | 2 | 0 |
| 0 | 0 | 1 | 1 | 1 | 2 | 2 | 2 | 0 | 0 | 0 | 1 |
| 0 | 0 | 0 | 0 | 0 | 0 | 0 | 0 | 1 | 1 | 1 | 1 |

| F13 | F14 | F15 | F16 | F17 | F18 | F19 | F20 | F21 | F22 | F23 | F24 |
|-----|-----|-----|-----|-----|-----|-----|-----|-----|-----|-----|-----|
| 1 | 2 | 0 | 1 | 1 | 2 | 0 | 1 | 2 | 0 | 1 | 2 |
| 1 | 1 | 2 | 2 | 0 | 0 | 1 | 1 | 1 | 2 | 2 | 2 |
| 1 | 1 | 1 | 1 | 2 | 2 | 2 | 2 | 2 | 2 | 2 | 2 |

**Contacts:**
(F1,F2), (F1,F3), (F1,F4), (F1,F9), (F1,F10), (F2,F4), (F2,F5), (F2,F10), (F2,F11), (F3,F4), (F3,F6), (F3,F9), (F3,F12), (F4,F5), (F4,F6), (F4,F7), (F4,F10), (F4,F12), (F4,F13), (F5,F7), (F5,F8), (F5,F11), (F5,F13), (F5,F14), (F6,F7), (F6,F12), (F6,F15), (F7,F8), (F7,F13), (F7,F15), (F7,F16), (F8,F14), (F8,F16), (F9,F10), (F9,F12), (F9,F17), (F9,F19), (F10,F11), (F10,F12), (F10,F13), (F10,F17), (F10,F18), (F10,F20), (F11,F13), (F11,F14), (F11,F18), (F11,F21), (F12,F13), (F12,F15), (F12,F19), (F12,F20), (F12,F22), (F13,F14), (F13,F15), (F13,F16), (F13,F20), (F13,F21), (F13,F23), (F14,F16), (F14,F21), (F14,F24), (F15,F16), (F15,F22), (F15,F23), (F16,F23), (F16,F24), (F17,F18), (F17,F19), (F17,F20), (F18,F20), (F18,F21), (F19,F20), (F19,F22), (F20,F21), (F20,F22), (F20,F23), (F21,F23), (F21,F24), (F22,F23), (F23,F24),

**Total:** $c(24) = y$ .

25 balls

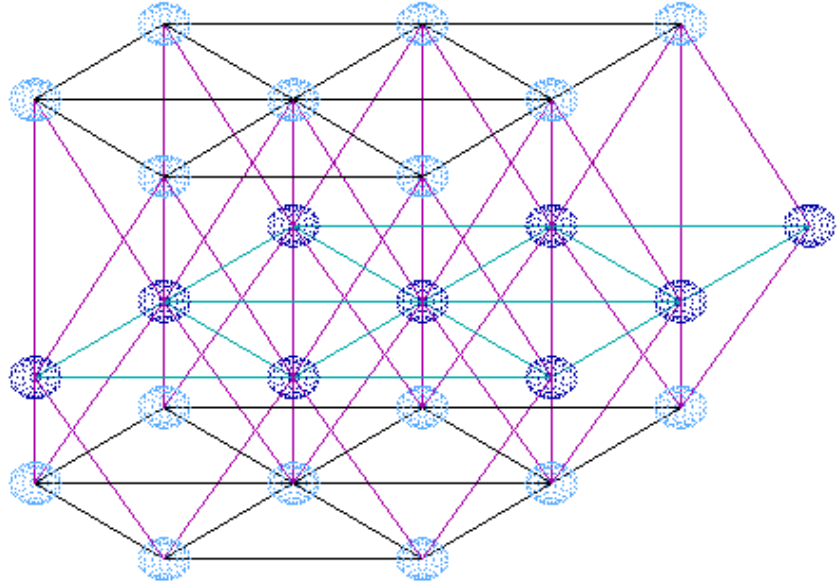

TEX25.acd

| F01 | F02 | F03 | F04 | F05 | F06 | F07 | F08 | F09 | F10 | F11 | F12 |
|---|---|---|---|---|---|---|---|---|---|---|---|
| 1 | 2 | 0 | 1 | 2 | 0 | 1 | 2 | 0 | 1 | 2 | 0 |
| 0 | 0 | 1 | 1 | 1 | 2 | 2 | 2 | 0 | 0 | 0 | 1 |
| 0 | 0 | 0 | 0 | 0 | 0 | 0 | 0 | 1 | 1 | 1 | 1 |

| F13 | F14 | F15 | F16 | F17 | F18 | F19 | F20 | F21 | F22 | F23 | F24 | F25 |
|---|---|---|---|---|---|---|---|---|---|---|---|---|
| 1 | 2 | 0 | 1 | 2 | 1 | 2 | 0 | 1 | 2 | 0 | 1 | 2 |
| 1 | 1 | 2 | 2 | 2 | 0 | 0 | 1 | 1 | 1 | 2 | 2 | 2 |
| 1 | 1 | 1 | 1 | 1 | 2 | 2 | 2 | 2 | 2 | 2 | 2 | 2 |

**Contacts:**
(F1,F2), (F1,F3), (F1,F4), (F1,F9), (F1,F10), (F2,F4), (F2,F5), (F2,F10), (F2,F11), (F3,F4), (F3,F6), (F3,F9), (F3,F12), (F4,F5), (F4,F6), (F4,F7), (F4,F10), (F4,F12), (F4,F13), (F5,F7), (F5,F8), (F5,F11), (F5,F13), (F5,F14), (F6,F7), (F6,F12), (F6,F15), (F7,F8), (F7,F13), (F7,F15), (F7,F16), (F8,F14), (F8,F16), (F8,F17), (F9,F10), (F9,F12), (F9,F18), (F9,F20), (F10,F11), (F10,F12), (F10,F13), (F10,F18), (F10,F19), (F10,F21), (F11,F13), (F11,F14), (F11,F19), (F11,F22), (F12,F13), (F12,F15), (F12,F20), (F12,F21), (F12,F23), (F13,F14), (F13,F15), (F13,F16), (F13,F21), (F13,F22), (F13,F24), (F14,F16), (F14,F17), (F14,F22), (F14,F25), (F15,F16), (F15,F23), (F15,F24), (F16,F17), (F16,F24), (F16,F25), (F17,F25), (F18,F19), (F18,F20), (F18,F21), (F19,F21), (F19,F22), (F20,F21), (F20,F23), (F21,F22), (F21,F23), (F21,F24), (F22,F24), (F22,F25), (F23,F24), (F24,F25),

**Total:** $c(25) = 84$.

26 balls

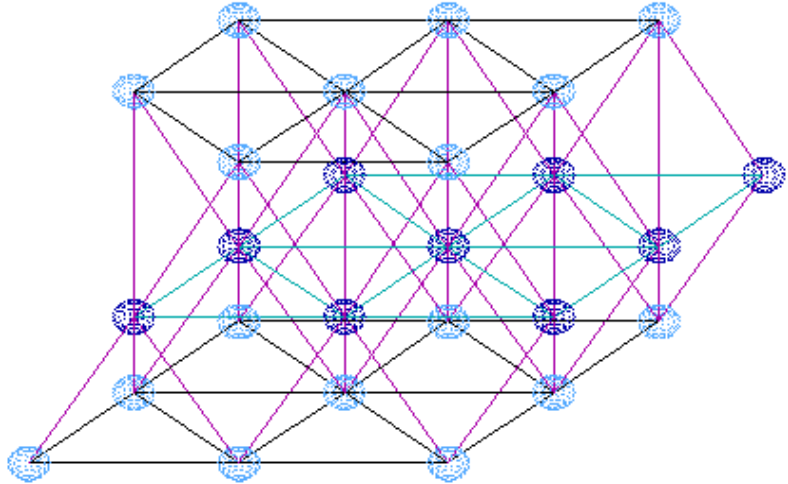

TEX26.acd

| F01 | F02 | F03 | F04 | F05 | F06 | F07 | F08 | F09 | F10 | F11 | F12 | F13 |
|---|---|---|---|---|---|---|---|---|---|---|---|---|
| 0 | 1 | 2 | 0 | 1 | 2 | 0 | 1 | 2 | 0 | 1 | 2 | 0 |
| 0 | 0 | 0 | 1 | 1 | 1 | 2 | 2 | 2 | 0 | 0 | 0 | 1 |
| 0 | 0 | 0 | 0 | 0 | 0 | 0 | 0 | 0 | 1 | 1 | 1 | 1 |

| F14 | F15 | F16 | F17 | F18 | F19 | F20 | F21 | F22 | F23 | F24 | F25 | F26 |
|---|---|---|---|---|---|---|---|---|---|---|---|---|
| 1 | 2 | 0 | 1 | 2 | 1 | 2 | 0 | 1 | 2 | 0 | 1 | 2 |
| 1 | 1 | 2 | 2 | 2 | 0 | 0 | 1 | 1 | 1 | 2 | 2 | 2 |
| 1 | 1 | 1 | 1 | 1 | 2 | 2 | 2 | 2 | 2 | 2 | 2 | 2 |

**Contacts:**
(F1,F2), (F1,F4), (F1,F10), (F2,F3), (F2,F4), (F2,F5), (F2,F10), (F2,F11), (F3,F5), (F3,F6), (F3,F11), (F3,F12), (F4,F5), (F4,F7), (F4,F10), (F4,F13), (F5,F6), (F5,F7), (F5,F8), (F5,F11), (F5,F13), (F5,F14), (F6,F8), (F6,F9), (F6,F12), (F6,F14), (F6,F15), (F7,F8), (F7,F13), (F7,F16), (F8,F9), (F8,F14), (F8,F16), (F8,F17), (F9,F15), (F9,F17), (F9,F18), (F10,F11), (F10,F13), (F10,F19), (F10,F21), (F11,F12), (F11,F13), (F11,F14), (F11,F19), (F11,F20), (F11,F22), (F12,F14), (F12,F15), (F12,F20), (F12,F23), (F13,F14), (F13,F16), (F13,F21), (F13,F22), (F13,F24), (F14,F15), (F14,F16), (F14,F17), (F14,F22), (F14,F23), (F14,F25), (F15,F17), (F15,F18), (F15,F23), (F15,F26), (F16,F17), (F16,F24), (F16,F25), (F17,F18), (F17,F25), (F17,F26), (F18,F26), (F19,F20), (F19,F21), (F19,F22), (F20,F22), (F20,F23), (F21,F22), (F21,F24), (F22,F23), (F22,F24), (F22,F25), (F23,F25), (F23,F26), (F24,F25), (F25,F26),

**Total:** $c(26) = \mathbf{87}$.

27 balls

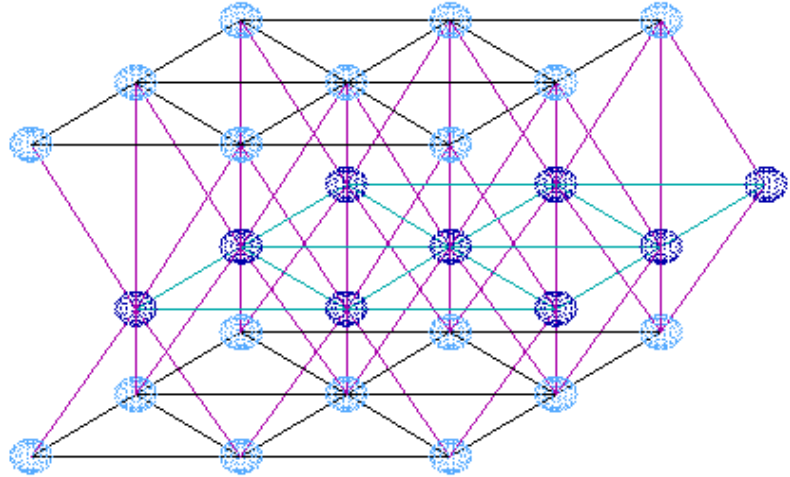

TEX27.acd

| F01 | F02 | F03 | F04 | F05 | F06 | F07 | F08 | F09 | F10 | F11 | F12 | F13 |
|---|---|---|---|---|---|---|---|---|---|---|---|---|
| 0 | 1 | 2 | 0 | 1 | 2 | 0 | 1 | 2 | 0 | 1 | 2 | 0 |
| 0 | 0 | 0 | 1 | 1 | 1 | 2 | 2 | 2 | 0 | 0 | 0 | 1 |
| 0 | 0 | 0 | 0 | 0 | 0 | 0 | 0 | 0 | 1 | 1 | 1 | 1 |

| F14 | F15 | F16 | F17 | F18 | F19 | F20 | F21 | F22 | F23 | F24 | F25 | F26 | F27 |
|---|---|---|---|---|---|---|---|---|---|---|---|---|---|
| 1 | 2 | 0 | 1 | 2 | 0 | 1 | 2 | 0 | 1 | 2 | 0 | 1 | 2 |
| 1 | 1 | 2 | 2 | 2 | 0 | 0 | 0 | 1 | 1 | 1 | 2 | 2 | 2 |
| 1 | 1 | 1 | 1 | 1 | 2 | 2 | 2 | 2 | 2 | 2 | 2 | 2 | 2 |

**Contacts:**
(F1,F2), (F1,F4), (F1,F10), (F2,F3), (F2,F4), (F2,F5), (F2,F10), (F2,F11), (F3,F5), (F3,F6), (F3,F11), (F3,F12), (F4,F5), (F4,F7), (F4,F10), (F4,F13), (F5,F6), (F5,F7), (F5,F8), (F5,F11), (F5,F13), (F5,F14), (F6,F8), (F6,F9), (F6,F12), (F6,F14), (F6,F15), (F7,F8), (F7,F13), (F7,F16), (F8,F9), (F8,F14), (F8,F16), (F8,F17), (F9,F15), (F9,F17), (F9,F18), (F10,F11), (F10,F13), (F10,F19), (F10,F20), (F10,F22), (F11,F12), (F11,F13), (F11,F14), (F11,F20), (F11,F21), (F11,F23), (F12,F14), (F12,F15), (F12,F21), (F12,F24), (F13,F14), (F13,F16), (F13,F22), (F13,F23), (F13,F25), (F14,F15), (F14,F16), (F14,F17), (F14,F23), (F14,F24), (F14,F26), (F15,F17), (F15,F18), (F15,F24), (F15,F27), (F16,F17), (F16,F25), (F16,F26), (F17,F18), (F17,F26), (F17,F27), (F18,F27), (F19,F20), (F19,F22), (F20,F21), (F20,F22), (F20,F23), (F21,F23), (F21,F24), (F22,F23), (F22,F25), (F23,F24), (F23,F25), (F23,F26), (F24,F26), (F24,F27), (F25,F26), (F26,F27),

**Total:** $c(27) = 90$.

# References

[1] **Arkus,N., Manoharan,V., Brenner,M.:** *Deriving Finite Sphere Packings*, SIAM J. Discrete Math. Vol. 25, No. 4, pp. 1860-1901, http://arxiv.org/abs/1011.5412v2 .

[2] **Bezdek,K.:** *Lectures on Sphere Arrangements - the Discrete Geometric Side*, Springer, 2013.

[3] **Bezdek,K., Reid,S.:** *Contact graphs of unit sphere packings revisited*, J. Geom. 104 (2013), no. 1, 57-83.

[4] **Bezdek,K., Szalkai,B., Szalkai,I.:** *On contact numbers of totally separable unit sphere packings*, Discrete Mathematics, Vol. 339, No. 2, (2015), Pages 668-676.

[5] **Bezdek,K., Khan,M.A.:** *Contact numbers for sphere packings*, http://arxiv.org/abs/1601.00145 .

[6] **Blair,D., Santangelo,C.D., Machta,J.:** *Packing Squares in a Torus*, J. Statistical Mechanics: Th. and Experiment, P01018 (2012), http://arxiv.org/abs/1110.5348, http://iopscience.iop.org/1742-5468/2012/01/P01018/article .

[7] **Harborth, H.:** *Lösung zu Problem* 664A, Elem. Math. 29, 14□15 (1974).

[8] **Holmes-Cerfon,M.:** *Enumerating nonlinearly rigid sphere packings*, http://arxiv.org/abs/1407.3285v2.

[9] **Hoy,R.S., Harwayne-Gidansky,J., O'Hern,C.S.:** *Structure of finite sphere packings via exact enumeration: Implications for colloidal crystal nucleation*, Physical Review E 85 (2012), 051403.

[10] **Neves, E.J.:** *A discrete variational problem related to Ising droplets at low temperatures*, J. Statistical Physics, 80 (1995), 103-123.

[11] **Reid, S.:** *Regular Totally Separable Sphere Packings*, RESEARCH June 2015, DOI: 10.13140/RG.2.1.2771.2161, http://arxiv.org/abs/1506.04171.

[12] **Reid, S.:** *On Contact Numbers of Finite Lattice Sphere Packings and the Maximal Coordination of Monatomic Crystals,* http://arxiv.org/pdf/1602.04246.pdf.

[13] **Szalkai,I.:** *Handling Multicomponent Systems in* $R^n$, J. Math.Chem. 25 (1999), 31-46.